\documentclass[11pt]{article}
\usepackage{graphicx}
\usepackage{amssymb,latexsym}
\usepackage{epstopdf}
\newtheorem{lemma}{Lemma}[section]
\newtheorem{theorem}[lemma]{Theorem}
\newtheorem{conj}[lemma]{Conjecture}
\newtheorem{prop}[lemma]{Proposition}

\newtheorem{cor}[lemma]{Corollary}
\newtheorem{remark}[lemma]{Remark}

\newcommand{\pf}{\noindent{\em Proof: }}
\newcommand{\epf}{\hfill\hbox{\rule{3pt}{6pt}}\\}

\newcommand{\forme}[1]{}
\begin{document}

\title{On Distance-Regular Graphs with Smallest Eigenvalue at Least $-m$ }

\author{
{\bf J.~H.~Koolen}\\
Pohang Mathematics Institute and Department of Mathematics,\\
 POSTECH\\
Hyoja-dong, Namgu,
Pohang 790-784 Korea\\
e-mail: koolen@postech.ac.kr\\
and\\
{\bf S.~Bang }
 \\
Department of Mathematics, Pusan National University\\
Geumjeong Gu, Busan 609-735 Korea\\
e-mail: sjbang3@pusan.ac.kr}

\date{\today}

\maketitle

\begin{abstract}
A non-complete  geometric distance-regular graph is the point
graph of a partial geometry in which the set of lines is a
set of Delsarte cliques. In this paper, we prove that for fixed
integer $m\geq 2$, there are only finitely many non-geometric
distance-regular graphs with smallest eigenvalue at least $-m$,
diameter at least three and intersection number $c_2 \geq 2$.
\end{abstract}

\section{Introduction}\label{intro}
In this paper, we will show that for fixed integer $m \geq 2$,
there are only finitely many distance-regular graphs with smallest
eigenvalue at least $-m$, diameter at least three, intersection
number $c_2$ at least two which are not the point graph of a
partial geometry. This result generalizes earlier results of
R. C. Bose \cite{bose} and A. Neumaier \cite{neumaier-79} for
strongly regular graphs, and of C. D. Godsil \cite{godsil-93} for
antipodal distance-regular graphs with diameter 3 (For
definitions, see next section.).

Note that any connected graph
has smallest eigenvalue at most $-1$ with equality if and only if the graph is complete. For connected regular graphs with smallest eigenvalue at least $-2$, it was shown by P. J. Cameron et al. \cite{cameron et al}, cf. \cite[Theorem 3.12.2]{bcn}, that either it is a line graph, a cocktail party graph or the number of vertices is at most 28.

To introduce the results of R. C. Bose, A. Neumaier and C. D. Godsil mentioned above, we will first introduce the notion of a Delsarte clique in a distance-regular graph.
Recall that a {\em clique} in a graph is a set of pairwise adjacent vertices.
 Let $\Gamma$ be a distance-regular graph with valency $k$, diameter $D \geq 2$ and the smallest eigenvalue $\theta_D$. Then any clique $C$ in $\Gamma$ contains at most $1 + {\displaystyle \frac{k}{-\theta_D}}$ vertices. This was shown by P. Delsarte \cite{Delsarte} for strongly regular graphs and
 C. D. Godsil generalized it to distance-regular graphs. A clique $C$ in $\Gamma$ is called a {\em Delsarte clique} if $C$ contains exactly $1 + {\displaystyle \frac{k}{-\theta_D}}$ vertices.
It is known that for any clique $C$ in $\Gamma$, the clique is Delsarte if and only if it is a completely regular code with covering radius $D-1$. Moreover, the outer distribution numbers of a Delsarte clique are completely determined by the intersection numbers of $\Gamma$ and hence do not depend on the specific Delsarte clique, cf. \cite[Section 13.7]{godsil-combin}. Note that for a distance-regular graph which contains a Delsarte clique, its smallest eigenvalue must be integral as ${\displaystyle \frac{k}{-\theta_D}}$ is integral.

C. D. Godsil \cite{godsil-93} introduced the following notion of a geometric distance-regular graph.
A non-complete distance-regular graph $\Gamma$ is called {\em geometric} if there exists a set of Delsarte cliques ${\mathcal C}$ such that each edge of $\Gamma$ lies in a unique $C \in {\mathcal C}$.  We will also say that  $\Gamma$ is {\em geometric with respect to  $\mathcal C$} in this case.

Examples of geometric distance-regular graphs include the Hamming
graphs, Johnson graphs, Grassmann graphs, dual polar graphs,
bilinear forms graphs and so on (See \cite[Chapter 9]{bcn} for
more information on these examples).

Note that a set $\mathcal C$ of Delsarte cliques for a geometric
distance-regular graph does not have to be unique. For example, in
the Johnson graph $J(2t,t), t \geq 2$, there are (exactly) two
different sets of Delsarte cliques such that $J(2t,t)$ is
geometric with respect to either one of them.

The definition of geometric distance-regular graphs for diameter
two is equivalent to the notion of geometric strongly regular
graphs as was introduced by R. C. Bose \cite{bose}. In the last
section, we will give more details.\\

Let $\Gamma$ be a geometric distance-regular graph with valency
$k$, diameter $D$ and the smallest eigenvalue $\theta_D$ and
assume that $\Gamma$ is geometric  with respect to $\mathcal C$.
Then the pair ${\mathcal G}:= (V(\Gamma), {\mathcal C})$ is a
partial geometry \footnote{A {\em partial geometry of order
$(s,t)$} is an incidence structure of points and lines such that
each line has $s+1$ points, each point is on $t+1$ lines and any
two distinct lines meet in at most one point.}(where a vertex $x$
is incident with a clique $C$ if and only if $x \in C$) of order
$(s,t)$ where $s= \displaystyle{\frac{k}{-\theta_D}}$ and $t =
-\theta_D -1$. This partial geometry is an example of a
distance-regular geometry as defined in F. De Clerck, S. De
Winter, E. Kuijken and C. Tonesi \cite{clerck-06}. The incidence
graph of $\mathcal{G}$ is an example of a distance-semiregular
graph as introduced by  H. Suzuki \cite{suzuki-semiregular}.

%Suppose that $\Gamma$ is a geometric distance-regular graph with valency $k\geq 2$, diameter $D\geq 2$ and smallest eigenvalue $\theta_D$, and that ${\cal C}$ is a set of Delsarte cliques partitioning the edge set of $\Gamma$.\\
%Then the pair $(\Gamma, \mathcal{C})$ is a Delsarte pair with parameters $(k,s_{\mathcal{C}},n_{\mathcal{C}})=(k,-\frac{k}{\theta_{D}},1)$ (cf. \cite[Definition 1.1]{DCG}), in other words, there exist constants $s_{\mathcal{C}},n_{\mathcal{C}}$, which depends only on $\mathcal{C}$, such that each clique in $\mathcal{C}$ has exactly $1+s_{\mathcal{C}}$ vertices and each edge of $\Gamma$ lies in $n_{\mathcal{C}}$ cliques in $\mathcal{C}$.\\ The partial linear space $(V(\Gamma), {\cal C})$ (where a vertex $x$ is incident with a clique $C$ if and only if $x \in C$) is a distance-regular geometry as defined in F. De Clerck, S. De Winter, E. Kuijken and C. Tonesi(\cite{clerck-06}).\\
%Note that the smallest eigenvalue $\theta_D$ of $\Gamma$ is an integer by (ii) of Definition \ref{def-geo}. Examples of geometric distance-regular graphs are  the Hamming graphs (and more general the regular
%$2D$-gons), the Johnson graphs, the Grassmann graphs and the bilinear forms graphs (\cite[p.506]{DCG}). For definitions and notation, we refer to Section 2.
\vskip0.01cm
A graph $\Gamma$ is called {\em coconnected} when its complement graph
(i.e., the graph with vertex set $V(\Gamma)$ whose edges are all the non-edges of $\Gamma$) is connected.
Note that the only non-complete distance-regular graphs which are not coconnected are the complete multipartite graphs
$K_{t \times n}$ with $t, n \geq 2$ (cf. \cite[Lemma 1.1.7]{bcn}).

In \cite{neumaier-79}, A. Neumaier has shown the following result.
\begin{theorem}\label{neum} {\em (cf. \cite[Theorem 4.6]{cameronsrg})}\\
Fix an integer $m\geq 2$. Then, there are only finitely many
coconnected non-geometric distance-regular graphs with smallest eigenvalue at least $-m$ and diameter two.
\end{theorem}

The next result is a generalization of Theorem \ref{neum}  to any
diameter at least two.

\begin{theorem}\label{main-diameter given}
Fix integers $m\geq 2$ and $D\geq 2$. Then there are only
finitely many coconnected non-geometric distance-regular graphs with smallest eigenvalue at least $-m$ and diameter
$D$.
\end{theorem}

In the next result, we show that we can replace the condition of a
fixed diameter by a condition on the intersection number $c_2$.

\begin{theorem}\label{main}
Fix an integer $m\geq 2$. Then there are only finitely many non-complete
coconnected non-geometric distance-regular graphs with smallest eigenvalue at least $-m$, and intersection number $c_2$
at least  $2$.
\end{theorem}

On \cite[p.130]{bcn}, they asked whether any distance-regular
graph with valency $k\geq 3$ and diameter $D\geq 3$ always has an
integral eigenvalue $\theta\neq k$. Theorem \ref{main} gives a
partial answer for this problem since the smallest eigenvalue of
any geometric distance-regular graph is integral. \vskip0.01cm In
a follow-up paper \cite{banitoconj}, we show that for fixed $k$ at
least three, there are only finitely many distance-regular graphs
with valency $k$. This implies that the condition $c_2 \geq 2$ can
be replaced by the condition that the valency is at least three in
Theorem \ref{main}. Note that the odd polygons are non-geometric distance-regular
graphs with smallest eigenvalue $>-2$.
%We conjecture that we may remove the assumption on $c_2$ in Theorem \ref{main}.  This is a slightly weaker conjecture than the conjecture of E. Bannai and T. Ito (\cite[p.237]{banito}) which states that for given integer $k \geq 3$, there are only finitely many distance-regular graphs with valency $k$.
\vskip0.01cm A. Neumaier \cite{neumaier-79} also showed that
except for a finite number of graphs, all geometric strongly
regular graphs with a given smallest eigenvalue are either Latin
square graphs or Steiner graphs (cf. \cite[Theorem
4.6]{cameronsrg}, We will discuss them in more detail in the last
section.). In \cite{wilson-74}, R.~M.~Wilson showed that there are
super-exponentially many Steiner graphs with parameters
$(v,k,\lambda,\mu)=(v,3s,s+3,9)$ for $v=\frac{(s+1)(2s+3)}{3}$
where $s\equiv 0$ or $2$ (mod $3$) and $s \geq 6$.  There are
super-exponentially many Latin square graphs for certain parameter
sets, see \cite[p. 210]{cameronsrg}. This shows that the
above-mentioned result of A. Neumaier is the best we can hope for
the case of distance-regular graphs of diameter two. We will
discuss the situation for geometric distance-regular graphs with
diameter at least three in the last section.

The paper is organized as follows. In Section 2, we give the
definitions. In Section 3, we give some useful results which we will use in the proofs of the above two theorems. In Section 4, we give some properties of
geometric distance-regular graphs.  In Section 5, we show, using results of K. Metsch, that the distance-regular
graphs with fixed smallest eigenvalue and intersection number  $c_2$ small compared to $a_1$, are geometric. This result in combination with the results in Section 3 implies Theorem 1.2.  In Section 6, we study the
distance-regular Terwilliger graphs with fixed smallest eigenvalue and $c_2 \geq 2$, and give a proof of Theorem 1.3. In the last section we will discuss the geometric distance-regular graphs in more detail and we end the paper with three conjectures.

\section{Definitions}
All the graphs considered in this paper are finite, undirected and
simple (for unexplained terminology and more details, see for
example \cite{bcn}). \vskip0.01cm

Let $\Gamma = (V(\Gamma), E(\Gamma))$ be a graph.
We write $x\sim_{\Gamma} y$ or simply $x\sim y$
if two vertices $x$ and $y$ are adjacent in $\Gamma$. For a vertex $x$ of $\Gamma$, let $\Gamma(x) := \{ y \mid y \sim x\}$, i.e. the set of neighbours of $x$ and
the {\em
valency} of $x$, denoted by $k(x)$, is the number of neighbours
of $x$, $|\Gamma(x)|$. The {\it local graph} of a vertex $x$ is
the subgraph of $\Gamma$ induced by $\Gamma(x)$.

We say that $\Gamma$ is {\em regular
with valency $k$} or {\em $k$-regular} if $k(x) = k$ for all
vertices $x$ of $\Gamma$.  A $k$-regular graph $\Gamma$ on $v$-vertices is called a {\em strongly regular graph } with parameters $(v,k,\lambda,\mu)$ if there are constants $\lambda$ and $\mu$ such that for any two distinct vertices $x$ and $y$, the number of common neighbours of $x$ and $y$ equals $\lambda$ if $x \sim y$ and $\mu$ otherwise.

The {\em adjacency matrix} $A = A(\Gamma)$ of $\Gamma$ is the $(|V(\Gamma)| \times |V(\Gamma)|)$-matrix whose  rows and the columns are indexed by
$V(\Gamma)$, and the $(x,y)$-entry of $A$
equals $1$ whenever $x \sim y$ and
$0$ otherwise. The eigenvalues of $\Gamma$ are the eigenvalues of $A$ and real as $A$ is a real symmetric matrix.

For the rest of this section, let $\Gamma$ be
a connected graph. The distance $d(x,y)$ between any two vertices $x,y$ of $\Gamma$
is the length of a shortest path between $x$ and $y$ in $\Gamma$.
The diameter of $\Gamma$ is the maximal distance occurring in
$\Gamma$ and we will denote this by $D = D(\Gamma)$.
%A connected regular graph with valency $k$ is called {\em amply regular} with parameters $(k,\lambda,\mu)$ if any two adjacent vertices have precisely $\lambda$ common neighbours and any two vertices of distance two have exactly $\mu$ common neighbours. Note that amply regular graphs with diameter two are exactly the same as connected non-complete strongly regular graphs.
%\vskip0.01cm
%{\bf Use uniformly regular partition instead of equitable partition : Delete $\rightarrow$ this part}A partition $\Pi= \{P_1, P_2, \ldots ,
%P_{\ell}\}$ of vertex set $V(\Gamma)$ is called {\em
%equitable} if there are constants $\beta_{ij}$ such that each vertex $x \in P_i$ has exactly $\beta_{ij}$ neighbours in $P_j$ ($1\leq i, j \leq \ell$).{\bf Delete $\leftarrow$}\\
For a vertex $x \in V(\Gamma)$, define $\Gamma_i(x)$ to be the set
of vertices which are at distance $i$ from $x~(0\le i\le D)$, % where
%$D:=\max\{d(x,y)\mid x,y\in V(\Gamma)\}$ is the diameter of
%$\Gamma$.
In addition, define $\Gamma_{-1}(x):=\emptyset$ and
$\Gamma_{D+1}(x) := \emptyset$. Note that $\Gamma(x)$ defined above is exactly the same as $
\Gamma_1(x)$ and

A connected graph $\Gamma$ with diameter $D$ is called {\em
antipodal} if for any vertices $x,y,z$ with $d(x,y)=D=d(y,z)$,
either $d(x,z)=D$ or $x=z$ hold.

A connected graph $\Gamma$ with diameter $D$ is called {\em
distance-regular} if there are integers $b_i,c_i$ $(0 \le i \leq
D)$ such that for any two vertices $x,y \in V(\Gamma)$ with
$d(x,y)=i$, there are precisely $c_i$ neighbours of $y$ in
$\Gamma_{i-1}(x)$ and $b_i$ neighbours of $y$ in $\Gamma_{i+1}(x)$
(cf. \cite[p.126]{bcn}). In particular, distance-regular graph
$\Gamma$ is regular with valency
$k := b_0$ and we define $a_i:=k-b_i-c_i$ for notational convenience. %Note that $a_i=|\Gamma(y)\cap \Gamma_i(x)|$ holds for any two vertices $x,y$ with $d(x,y)=i$ ($0\leq i\leq D$).
The numbers $a_i$, $b_{i}$ and $c_i~(0\leq i\leq D)$ are called the {\em
intersection numbers} of $\Gamma$. Note that $b_D=c_0=a_0=0$ and $c_1=1$.
The intersection numbers of a distance-regular graph $\Gamma$ with diameter $D$ and valency $k$ satisfy
(cf. \cite[Proposition 4.1.6]{bcn})\\
(i) $k=b_0> b_1\geq \cdots \geq b_{D-1}$;\\
(ii) $1=c_1\leq c_2\leq \cdots \leq c_{D}$;\\
(iii) $b_i\ge c_j$ \mbox{ if }$i+j\le D$.\\
Moreover, if we fix a vertex $x$ of $\Gamma$ and define $k_i := |\Gamma_i(x)|$, then $k_i$ does not depend on the choice of $x$ as $c_{i+1} k_{i+1} =
b_i k_i$ hold for $i =1, 2, \ldots D-1$.
Note that the non-complete connected strongly regular graphs are exactly the distance-regular graphs with diameter two.

For the rest of this section, let $\Gamma$ be a distance-regular graph with diameter $D$. Let $C,C' \subseteq V(\Gamma)$ be non-empty subsets and $x$ be a vertex of $\Gamma$. We write  $d(x,C):=\min\{d(x,y)\mid y\in C\}$ and $d(C',C):=\min\{d(x,y)\mid x\in C',\,y\in C\}$.
The {\em covering radius} of $C$, denoted by $\rho(C)$ is defined as $\rho(C):=\max\{ d(x, C) \mid x \in V(\Gamma)\}$, and define $C_i :=\{x\in V(\Gamma) \mid d(x,C)=i\}~(0\le i \le \rho(C))$. For $x$ a vertex of $\Gamma$ and $C$ a non-empty subset of $V(\Gamma)$, we write
$B_{xi}(C):=\left| C\cap \Gamma_i(x)\right|$.
The numbers $B_{xi}(C)$, $ i=0, 1,\ldots, D $, are called the {\em outer distribution numbers} of $C$.

 A non-empty subset $C\subseteq V(\Gamma)$ with covering radius $\rho$, is called a {\em completely regular code}, if the outer distribution number $B_{xi}(C)$ only depends on $i$ and $d(x,C)$, that is, there exist numbers $e_{\ell\,i}$  $(\ell = 0, 1, \ldots, \rho, \ i= 0, 1, \ldots, D)$ such that for all vertices $x$ of $\Gamma$ and $i \in \{0, 1, \ldots D\}$, we have $B_{xi}(C) = e_{\ell\,i}$ where  $\ell = d(x,C)$. We refer to the numbers $e_{\ell\,i}$ as the {\em outer distribution numbers} of the completely regular code $C$ and write $\psi_i(C) := e_{ii}$ for $i =0,1,\ldots, \rho(C)$.

 A partition $\Pi= \{P_1, P_2, \ldots ,
P_{\ell}\}$ of $V(\Gamma)$ into non-empty parts, is called {\em
equitable} if there are constants $\beta_{ij}$ such that each vertex $x \in P_i$ has exactly $\beta_{ij}$ neighbours in $P_j$ ($1\leq i, j \leq \ell$).
The {\em quotient matrix} of $\Pi$ is the $(\ell \times \ell)$-matrix  $Q= Q(\Pi)$ defined by $Q_{ij} := \beta_{ij}$ for $1 \leq i,j \leq \ell$.
An equitable partition $\Pi$ of $V(\Gamma)$ is called a {\em uniformly regular partition} if there exist numbers $e_{01}$ and $e_{11}$ such for any $C\in \Pi$ and $x\in V(\Gamma)$, the number $ B_{x1}(C) $ is equal to $ e_{\ell1}$ if $\ell \in \{0, 1\}$ and zero otherwise, where $\ell = d(x, C)$.
Moreover, a uniformly regular partition $\Pi$ of $V(\Gamma)$ is called a {\em completely regular partition} if each $C\in \Pi$ is a completely regular code and the outer distribution numbers for all $C \in \Pi$ are the same.
\vskip0.01cm

Suppose that $\Gamma$ is a distance-regular graph with valency $k\ge 2$ and
diameter $D\ge 2$, and $A=A(\Gamma)$ is the adjacency matrix
of $\Gamma$. It is well-known that $\Gamma$ has exactly $D+1$ distinct eigenvalues, $k=\theta_0>\theta_1>\cdots>\theta_D$ (\cite[p.128]{bcn}). For an eigenvalue $\theta$ of $\Gamma$, the sequence $u_i=u_i(\theta)$ $(0\leq i\leq D)$ satisfying
\begin{equation}\label{rec-d}
u_0=1,~~u_1=\frac{\theta}{k},~~c_i u_{i-1} + a_i u_i + b_i u_{i+1} = \theta u_i
\end{equation}
is called the {\em standard sequence} corresponding to the eigenvalue $\theta$ (\cite[p.128]{bcn}). It is known that the standard sequence corresponding to $\theta_i$ has exactly $i$ sign changes (\cite[Corollary 4.1.2]{bcn}).
\vskip0.01cm

% Let $e, n,D$ be integers with $e\geq 2$ and $n\geq 2D\geq 4$. The Grassmann graph $J_q(n,D)$ is the graph whose vertices are the subspaces of dimension $D$ of a given vector space $V$ with dimension $n$ over $\mathbb{F}_q$ the finite field on $q$-elements, where two vertices are adjacent if they meet in a subspace of $V$ of dimension $D-1$. Note that $J_q(n,D)$ is a distance-regular graph with valency \[\frac{q (q^D-1) (q^{n-D}-1)}{(q-1)^2}.\] The bilinear forms graph $H_q(n,D)$ is the induced subgraph of $J_q(n+D,D)$ induced on the vertices which meet non-trivially a fixed subspace $N$ of dimension $n$. The Johnson graph $J(n,D)$ is the graph with the vertex set as the $D$-subsets of given $n$-set $X$ whose two vertices are adjacent if they meet in exactly $(D-1)$-elements. The folded Johnson graph $\bar{J}(4D,2D)$ is the graph with the vertex set as the set of partitions of $4D$-set $X$ into two $2D$-sets, where two partitions are adjacent whenever their common refinement is a partition of $X$ into four sets of sizes $1,1,2D-1,2D-1$. The Hamming graph $H(e,D)$ is the graph whose vertices are the cartesian product of $D$-copies of $e$-set $X$, where two vertices are adjacent if they differ in precisely one coordinate. Note that $H_q(n,D)$, $J(n,D)$, $\bar{J}(4D,2D)$ and $H(e,D)$ are distance-regular graphs with valencies
%\[\frac{(q^D-1) (q^n-1)}{q-1}, ~n(n-D),~4D^2\mbox{ and }D(e-1),\] respectively (For more information, see \cite[Chapter 9]{bcn}.).

\section{Some Useful Results}

In this section, we list a number of results which will be used in this paper.

First, we show that for a connected $k$-regular graph
of diameter at least three, it has an eigenvalue $\theta \neq k$ satisfying $|\theta| > \sqrt{\frac{k}{2}}$.

\begin{lemma}\label{theta_1}
 Suppose that $\Gamma$ is a connected regular graph with valency $k\ge 2$ and diameter $D\ge 3$. Then $\Gamma$ has an eigenvalue $\theta$ different from $k$ satisfying $$|\theta|> \sqrt{\frac{k}{2}}.$$
\end{lemma}

\pf Suppose that $\Gamma$ has $n$ vertices and let $k= \eta_1 > \eta_2 \geq \ldots \geq \eta_n$ be the eigenvalues of $\Gamma$. Let $A$ be the adjacency matrix of $\Gamma$.
Then it is well-known (cf. \cite[Lemma 2.5]{biggs}) that
$$\sum_{i=1}^n \eta_i^2 = \mbox{Tr}(A^2) = nk,$$ where Tr$(A^2)$ is the trace of $A^2$.
Let $x,y$ be two vertices in $\Gamma $ at distance 3.
Then $n \geq
|\Gamma(x) \cup \{x\} \cup \Gamma(y) \cup \{y\}|=2k +2$.
It follows that $\sum_{i=2}^n \eta_i^2 = (n-k)k > \frac{n}{2} k$. This immediately implies the lemma.
\epf

Next we show two results on distance-regular graphs with smallest eigenvalue at least $-m$.

\begin{lemma}\label{lowerbound}
Fix an integer $m\geq 2$. Suppose that $\Gamma$ is a distance-regular
graph with smallest eigenvalue at least $-m$, valency $k\ge 2$, diameter $D\ge 2$ and intersection
number $a_1$. Then
$$%\label{a_1-lowerbound}
k< m(a_1+m).
$$
\end{lemma}
\pf Let $\theta_D$ be the smallest eigenvalue of $\Gamma$. Since the standard sequence of $\theta_D$ has exactly $D$ sign changes, it follows that $u_2(\theta_D)>0$. By (\ref{rec-d}), we have  $u_2(\theta_D)=\frac{\theta_D^2-a_1 \theta_D-k}{k(k-a_1-1)}$, and by using $\theta_D\geq -m$, %Inequality (\ref{a_1-lowerbound})
the lemma  follows.
\epf

\begin{theorem}\label{nonterw}
Fix a  real number $\epsilon$, $0< \epsilon <1$, and integers
$m\geq 2$ and $D\geq 3$. Then there are only finitely many
distance-regular graphs with smallest eigenvalue at
least $-m$, diameter $D$ and intersection numbers $a_1$, $c_2$ satisfying
$c_2\geq \epsilon a_1$.
\end{theorem}

\pf For a given real number $\epsilon$, $0< \epsilon <1$ and for given integers
$m\geq 2$ and $D\geq 3$, suppose that $\Gamma$ is a
distance-regular graph with valency $k \geq 2$ and diameter $D$ such that its smallest eigenvalue is at
least $ -m$ and its intersection numbers $a_1$, $c_2$ satisfy
$c_2\geq \epsilon a_1$.
We will show that \begin{equation}\label{k-equation}
k < D^2
\left(\frac{2m^2}{\epsilon}\right)^{2D+4}
\end{equation} holds, from which the theorem immediately follows, as the diameter $D$ is fixed.

By Lemma \ref{lowerbound},  we have $k < m(a_1 +m)$ and hence, if $a_1 = 0$, then $k < m^2$ and we are done.
Hence, in order to show Inequality (\ref{k-equation}), we may assume that   $a_1\ne 0$ and $k\ge2m^2\left(\frac{2m^2}{\epsilon}\right)^{2 D}D^2$.\\
As $k < m(a_1 +m)$ and $c_2 \geq  \epsilon a_1$, we find
\begin{equation} \label{b1c2}
\frac{b_1}{c_2} < \frac{(m-1)(a_1+m+1)}{
\epsilon a_1}<\frac{2m^2}{\epsilon}.
\end{equation}
As $b_{i-1} \geq b_{i}$ and $c_i \leq c_{i+1}$ hold for all $i=1, \ldots, D-1$ and using Inequality (\ref{b1c2}), we obtain
$$k_{i+1}= \frac{b_i}{c_{i+1}}k_i \leq \frac{b_1}{c_2} k_i < \left(\frac{2
m^2}{\epsilon}\right)^i k~~(i=1,\ldots,D-1). $$
This implies
\begin{equation}\label{vertexbound}
|V(\Gamma)|= \sum_{i=0}^{D}k_i < \left(\frac{2
m^2}{\epsilon}\right)^{D} D k.
\end{equation}

Now  Lemma \ref{theta_1} implies that $\theta_1>\sqrt{\frac{k}{2}}$ as
$k>2m^2$, and thus from (\ref{vertexbound}) we obtain
\begin{equation}\label{tr(A')}
0=\mbox{Tr}(A) = \sum_{i=0}^D m_i \theta_i >m_1\theta_1-m|V(\Gamma)|>\sqrt{\frac{k}{2}}m_1- m
\left(\frac{2m^2}{\epsilon}\right)^{D}Dk.
\end{equation}

By \cite[Proposition 4.4.8]{bcn} and by Inequality (\ref{tr(A')}) with $k\geq 2m^2\left(\frac{2m^2}{\epsilon}\right)^{2 D}D^2$, it follows $2<m_1<k$. This in turn implies, by \cite[Theorem 4.4.4]{bcn} that any local graph of
$\Gamma$ has an eigenvalue $\eta$, where $\eta:=-1-\frac{b_1}{\theta_1+1}$.
As any eigenvalue of a subgraph of $\Gamma$ is at least the smallest eigenvalue of $\Gamma$, $\eta \geq -m$ holds. As $b_1 \geq c_2$ holds by $D \geq 3$, it follows that
$k<m(a_1+m)\leq 2m^2 a_1\leq \frac{2m^2}{\epsilon}c_2\leq
\frac{2m^2}{\epsilon}b_1$ which in turn implies
\begin{eqnarray}
\theta_1&\ge &\frac{b_1}{m-1}-1>\frac{\epsilon k}{2m^2(m-1)}-1\ge
\frac{\epsilon k}{2m^3}.\label{case2-theta_1-bd}
\end{eqnarray}
Moreover, by Inequalities (\ref{vertexbound}), (\ref{tr(A')}) and (\ref{case2-theta_1-bd}),
$$m_1\frac{\epsilon k}{2m^3} < m_1 \theta_1 < m
\left(\frac{2m^2}{\epsilon}\right)^{D}Dk,$$
and this implies
$$m_1<Dm^2\left(\frac{2m^2}{\epsilon}\right)^{D+1}<D \left(\frac{2m^2}{\epsilon}\right)^{D+2}.$$
As $2k\le (m_1-1)(m_1+2)$, by  \cite[Theorem 5.3.2]{bcn},  we find
$2k\le (m_1-1)(m_1+2)< 2D^2
\left(\frac{2m^2}{\epsilon}\right)^{2D+4}$ holds, and this
completes the proof. \epf

Now, we consider diameter bounds for distance-regular graphs. A special case of the diameter bound given by A. A. Ivanov (\cite{aaivanov}, cf. \cite[Theorem 5.9.8]{bcn}) gives that the diameter of a distance-regular graph with valency $k$ and $c_2\geq 2$ satisfies $D\leq 4^k$.
As there are only finitely many connected non-isomorphic $k$-regular graphs with diameter at most $4^k$, we obtain:
\begin{theorem}\label{ivanov}
Fix an integer $k\ge 3$. Then there are only finitely many distance-regular graphs with valency $k$ and intersection number $c_2\ge 2$.
\end{theorem}

The next diameter bound is due to P. Terwilliger (\cite{terw-85}, cf. \cite[Corollary 5.2.2]{bcn}).

\begin{theorem}\label{diamterw}
 Suppose that $\Gamma$ is a distance-regular graph with valency $k\geq 2$, diameter $D\geq 2$ and intersection numbers $a_1$ and $c_D$. If $\Gamma$ contains an induced quadrangle, then
$$D \leq \frac{k+ c_D}{a_1 +2} \leq \frac{2k}{a_1+2}.$$
\end{theorem}

\begin{prop}\label{diam-m}
Fix an integer $m\geq 2$. Suppose that $\Gamma$ is a distance-regular
graph with smallest eigenvalue at least $-m$, valency $k\geq 2$, diameter $D\geq 2$ and intersection
number $a_1$. If $\Gamma$ contains an induced quadrangle, then
$$
D <\frac{2m(a_1+m)}{a_1+2}\le m^2.
$$
\end{prop}
\pf The result follows immediately from Lemma \ref{lowerbound} and Theorem \ref{diamterw}.
\epf

\section{Properties of Geometric Distance-Regular Graphs}
In this section, we will give some properties of  geometric distance-regular graphs,
which we will need later in this paper.\\
First, we recall the notion of a Delsarte pair as was introduced
by S. Bang, A. Hiraki and J. Koolen \cite{DCG}. Suppose that
$\Gamma$ is a non-complete distance-regular graph with valency
$k$, diameter $D$ and the smallest eigenvalue $\theta_D$, and that
$\mathcal C$ is a set of Delsarte cliques in $\Gamma$. Then the
pair $(\Gamma, {\mathcal C})$ is called a {\em  Delsarte pair with
parameters $(k, \frac{k}{-\theta_D},n_{\mathcal{C}})$}, if every
edge lies in exactly $n_{\mathcal{C}}$ cliques of $\mathcal C$,
where $n_{\mathcal{C}}\geq 1$ (cf. \cite[Definition 1.1]{DCG}).

Suppose that $\Gamma$ is a geometric distance-regular graph with respect to $\mathcal{C}$, that is,
 $(\Gamma, \mathcal{C})$ is a Delsarte pair with $n_{\mathcal{C}} = 1$.
Recall that for a Delsarte clique $C$, the numbers $\psi_i(C)$ for $i=1,2, \ldots, D-1$ do only depend
 on the intersection numbers of $\Gamma$.   Hence, we can write $\psi_i := \psi_i(C)$ for $C \in \mathcal{C}$.
  In \cite[Lemma 4.1,~Proposition 4.2 (i)]{DCG}, we showed that for any fixed  integer $j~,1\leq j\leq D$
  and for any vertices $x,y$ at distance $j$, the number of cliques $C \in \mathcal{C}$
that contain $y$ and satisfy $d(x,C)=j-1$ is dependent only on $j$
and $n_{\mathcal C}$ (and not on the particular pair of $x,y$ at
distance $j$) and we denote this number by
$\tau_j:=\tau_j(\mathcal{C})$. Note that the numbers $\psi_i$ and
$\tau_j$ do not depend on the particular set of Delsarte
cliques.\\
Therefore the next lemma is a direct consequence of
\cite[Proposition 4.2 (i)]{DCG} with $n_{\mathcal{C}}=1$.
\begin{lemma} \label{geometricpar}
Suppose that $\Gamma$ is a geometric distance-regular graph with
valency $k\geq 2$, diameter $D\geq 2$ and smallest eigenvalue
$\theta_D$. Then the following hold :
\begin{eqnarray*}
b_i & =& -(\theta_D + \tau_i)\left( 1+ \frac{k}{-\theta_D} - \psi_i \right)~~(1\leq i\leq D-1),\\
c_i &= &\tau_i \psi_{i-1}~~(1\leq i\leq D).
\end{eqnarray*}
%In particular, $c_D \geq -\theta_D > 1$ holds.
\end{lemma}

A distance-regular graph does not need to be geometric if its intersection numbers satisfy the equations in Lemma
\ref{geometricpar} for some non-negative integers $\tau_i, \psi_{i-1}~(i=1,2,
\ldots, D)$. For example, the Doob graph $D(s,c)$ (the direct
product of $s$ Shrikhande graphs and $c$ 4-cliques with $s \geq
1$) is non-geometric and it has the same intersection numbers with
the Hamming graph $H(2s+c,4)$, which is geometric.\\
Part (ii) of the next lemma is a direct consequence of \cite[Lemma 5.1
(iii)]{DCG}, but  we include its proof for completeness.

\begin{lemma}\label{tau}
Suppose that $\Gamma$ is a geometric distance-regular graph with
diameter $D\geq 2$ and intersection number $c_2$ at least two.
%(i) If $\psi_1 \geq 2$, then $\psi_{i} > \psi_{i-1}$ for all $2\leq i \leq D-1$.\\
Then
the following holds:\\
(i) $\tau_2 \geq \psi_1$.\\
(ii) $\Gamma$ contains an induced quadrangle.
\end{lemma}
\pf Let $\mathcal C$ be a set of Delsarte cliques with respect of
which $\Gamma$ is geometric. Let $x,y$ be any two vertices of
$\Gamma$ at distance 2. Let $z$ be a common neighbour of $x$ and
$y$, and let $C$ be the unique clique in $\mathcal{C}$ containing
$x$ and $z$. Put $\Gamma(y) \cap C:=\{z_1, \ldots, z_{\psi_1} \}$.
Then each edge $\{y,z_i\} $ lies in a unique clique $C_{(i)} \in
{\mathcal C}$, $1\leq i\leq \psi_1$.
\vskip0.01cm
(i): As $C_{(i)} \neq C_{(j)}$ holds for any $i \neq j$, the result (i) follows.\\
(ii): If $\psi_1=1$,  then the local graph of any vertex is
disjoint union of cliques, and the result holds as $c_2\geq 2$.
Now, suppose $\psi_1\geq 2$. Then there exists a vertex $w \in
C_{(\psi_1)}\cap \Gamma(x)$ such that $w \not\sim z_1$ holds, by
the definition of $\psi_1$ and $z_1\sim y$. Hence the subgraph
induced on $x,y,z_1,w$ is a quadrangle. This shows (ii). \epf

%Then, for any
%$z'\in \left(C\cap \Gamma(y)\right)\setminus \{z\}$, and the
%Delsarte clique $C' \in \mathcal{C}$ containing $y$ and $z'$, there exists a
%vertex $w$ in $C'\cap \Gamma(x)$ satisfying $w\not\sim z$ by
%$z\sim y$ and the definition of $\psi_1$. Hence, the set
%$\{x,y,z,w\}$ induces a quadrangle.

By applying Lemma \ref{tau} (ii) to Proposition \ref{diam-m}, we obtain the following diameter bound
for geometric distance-regular graphs.

\begin{prop}\label{geom-diam}
Fix an integer $m\geq 2$. Suppose that $\Gamma$ is a geometric
distance-regular graph with smallest eigenvalue $-m$, diameter $D\geq 2$ and intersection number $c_2$.
If $c_2\geq 2$ holds, then
$$
D <m^2.
$$
\end{prop}

\section{Distance-Regular Graphs with Small $c_2$}
In this section, we will show that any distance-regular graphs with intersection number $c_2$ much smaller then  $a_1$ are geometric. Also we will show Theorem 1.2.

To show any distance-regular graphs with intersection number $c_2$ much smaller then  $a_1$ are geometric, we first state the following result of K. Metsch.

\begin{prop}\label{result2.2} (\cite[Result 2.1]{metsch-99}){\ \\}
Let $k\geq 2,~\mu\ge 1,~\lambda\geq 0,~s\geq 1$ be integers.
Suppose that $\Gamma$ is a $k$-regular graph such that
any two non-adjacent vertices have at most $\mu$ common neighbours, and any two adjacent vertices have exactly $\lambda$ common neighbours. Define a line to be a maximal clique $C$ in $\Gamma$ such that $C$ has at least $\lambda+2-(s-1)(\mu-1)$ vertices. If\\
\noindent(i) $\lambda>(2s-1)(\mu-1)-1$ and\\
(ii) $k<(s+1)(\lambda+1)-\frac{1}{2}s(s+1)(\mu-1)$\\
all hold, then every vertex is in at most $s$ lines, and each edge lies in a unique line.
\end{prop}

\noindent As a consequence of this result, we obtain the following
result on distance-regular graphs.

\begin{prop}\label{amply-regular thm}
Fix an integer $m\geq 2$. Suppose that $\Gamma$ is a
distance-regular graph with smallest eigenvalue at least $-m$,
valency $k\geq 3$, diameter $D\geq 2$ and intersection numbers
$a_1$ and $c_2$. We define a line to be a maximal clique $C$ in
$\Gamma$ such that $C$ has at least $a_1+2-(m-1)(c_2-1)$ vertices.
If $a_1>m^2c_2$, then every vertex is in at most $m$ lines, and
each edge lies in a unique line.
\end{prop}

\pf Suppose that $\Gamma$ is a
distance-regular graph with smallest eigenvalue at least $-m$,
valency $k\geq 3$, diameter $D\geq 2$ satisfying
$a_1>m^2c_2$. By Proposition \ref{result2.2}, it is enough to show
that the following inequalities all hold:
\begin{eqnarray}
a_1&>&(2m-1)(c_2-1)-1, \label{2.1-a}\\
k&<&(m+1)(a_1+1)-\frac{1}{2}m(m+1)(c_2-1)\label{2.1-b}.
\end{eqnarray}
Inequality (\ref{2.1-a}) follows immediately from $a_1>m^2c_2$. To see that  Inequality (\ref{2.1-b}) holds, note that
$$
m(a_1 + m) \leq (m+1)(a_1+1)-\frac{1}{2}m(m+1)(c_2-1) $$ holds,
from which Inequality (\ref{2.1-b}) follows by Lemma
\ref{lowerbound} \epf

The following theorem is the main result of this section.

\begin{theorem}\label{geometric}
Fix an integer $m\geq 2$. Suppose that $\Gamma$ is a
distance-regular graph with diameter $D\ge 2$ such that its
smallest eigenvalue $\theta_D$ satisfies $-m\leq \theta_D< 1-m$.
If intersection numbers $a_1$ and $c_2$ satisfy $a_1>m^2c_2$, then
$\Gamma$ is geometric (and $ \theta_D=-m$).
\end{theorem}

\pf Suppose that $\Gamma$ is a
distance-regular graph with diameter $D\ge 2$ such that its
smallest eigenvalue $\theta_D$ satisfies $-m\leq \theta_D< 1-m$, satisfying
$a_1>m^2c_2$. Define lines as in Proposition \ref{amply-regular
thm}, and for each vertex $x\in V(\Gamma)$, let $M_x$ be the
number of lines containing $x$.  By Proposition \ref{amply-regular
thm}, we have $M_x \leq m$. On the other hand, as any maximal
clique $C$ satisfies $|C| \leq 1 + \frac{k}{-\theta_D} < 1 +
\frac{k}{m-1}$,
$$M_x > \frac{k}{\left(1+ \frac{k}{m-1}\right) -1} = m-1$$
follows as every edge lies in a unique line, by Proposition
\ref{amply-regular thm}. Hence, $M_x = m$ for all $x$.

Let $B$  be the vertex-line incidence matrix (i.e.  the
$(0,1)$-matrix whose rows and columns are indexed by the vertex
set and the set of lines of $\Gamma$, respectively where
$(x,C)$-entry of $B$ is $1$ if the vertex $x$ is contained in the
line $C$ and $0$ otherwise). Then $B B^T=A(\Gamma) + mI$ holds,
where $B^T$ is the transpose of $B$ and $I$ is the identity
matrix. Since each line contains more than $m$ vertices, by $ a_1
+ 2 - (m-1)(c_2-1) > m c_2 \geq  m$, the matrix $BB^T$ is
singular. This implies that $0$ is an eigenvalue of $BB^T$ and
thus $-m$ is an eigenvalue of $A$. As  $\theta_D\geq -m$,
$\theta_D = -m$. But then, every line has exactly $1 +\frac{k}{m}$
vertices as any maximal clique has cardinality at most $1 +
\frac{k}{-\theta_D} = 1 + \frac{k}{m}$, $M_x =m$ for all $x$ and
each vertex lies in a unique line. This proves that $\Gamma$ is
geometric with $\theta_D=-m$. \epf

%\section{Distance-Regular Graphs with Large $c_2$}

%In this section we show  that for given integers $m\geq 2$ and
%$D\geq 3$, there are only finitely many distance-regular graphs
%with diameter $D$, the smallest eigenvalue at least $-m$ and
%intersection number $c_2$ which is large in comparison with $a_1$.
%Moreover, we will use this result to prove Theorem \ref{main-diameter given}.\\

%Let $\Gamma$ denote a regular graph with valency $k\geq 2$ and $k=\theta_0>\theta_1>\cdots>\theta_n$ be the eigenvalues of $\Gamma$.
%A sequence of vertices $W=v_0,v_1, \cdots , v_{\ell}$ with $\ell\geq 1$ an integer, which are not necessarily mutually distinct, is called a {\it $\ell$-walk} if $v_i$ and $v_{i+1}$ are
%adjacent for all $i=0,\cdots , {\ell}-1$. If $v_0=v_{\ell}$ then $W$ is
%called a {\it closed $\ell$-walk}. It is well-known that for each $1\leq j\leq n$,
%\begin{equation}\label{tr}
%\sum_{i=0}^{n}m_i\theta_i^j=\mbox{Tr}(A^j)= \left| \{ W \mid W \mbox{ is a
%closed $j$-walk } \}\right |
%\end{equation}
%follows (cf. \cite[Lemma 2.5]{biggs}), where $\mbox{Tr}(A^j)$ is the trace of $A^j$ and $m_i$ is the multiplicity of $\theta_i~(0\leq i\leq n)$.

\noindent {\em Proof of Theorem \ref{main-diameter given}:} For
given integers $m\geq 2$ and $D\geq 2$, let $\Gamma$ be a
non-geometric distance-regular graph with diameter $D$ and smallest eigenvalue at least $-m$.
By Theorem \ref{geometric}, the
intersection numbers $a_1$ and $c_2$ of $\Gamma$ satisfy $c_2 \geq
\frac{1}{m^2}a_1$. The result now immediately follows from Theorem
\ref{nonterw} and Theorem \ref{neum}. \epf

\section{Terwilliger Graphs and Proof of Theorem 1.3}
In this section, we will show Theorem 1.3. Before we do this, we
first need to consider Terwilliger graphs which are
distance-regular (i.e., distance-regular graphs without induced
quadrangles). In \cite[p.36, Problems (ii)]{bcn}, it was asked whether it is possible to classify
the distance-regular Terwilliger graphs with
intersection number $c_2\ge 2$. We will show in Theorem
\ref{terwthm} that for a fixed integer $m\geq 2$, there are only
finitely many distance-regular Terwilliger graphs with $c_2 \geq
2$ and smallest eigenvalue at least $-m$.

For an integer $\alpha\ge 1$, a {\em $\alpha$-clique extension} of a graph $\Gamma$ is the graph
$\Sigma$ obtained from $\Gamma$ by replacing each vertex $x\in V(\Gamma)$ by a clique $C(x)$ of $\alpha$ vertices,
where for any $x,y \in V(\Gamma)$, $x' \in C(x) $ and $y' \in C(y)$, $x \sim_{\Gamma} y$ if and only if
$x' \sim_{\Sigma} y'$.

\begin{prop}\label{a-clique-ev}
 Suppose that $\Gamma$ is a graph and that $\Sigma$ is the $\alpha$-clique extension of $\Gamma$ for an integer
 $\alpha \geq 1$. Then for each eigenvalue $\theta$ of $\Gamma$, $(\alpha (\theta+1)-1)$ is an eigenvalue of $\Sigma$.
\end{prop}

\pf Let $\Pi$ be the set of all $\alpha$-cliques of $\Sigma$ which correspond to the vertices of $\Gamma$, respectively.
Then $\Pi$ is a uniformly regular partition of $V(\Sigma)$ and the quotient matrix $Q(\Pi)$ of $\Pi$ satisfies
$Q(\Pi)=\alpha A(\Gamma)+(\alpha-1)I$. Hence if $\theta$ is an eigenvalue of $\Gamma$, then $(\alpha(\theta+1)-1)$ is an
eigenvalue of $Q(\Pi)$. By \cite[Lemma 5.2.2]{godsil-combin}, $(\alpha (\theta+1)-1)$ is also an eigenvalue of $\Sigma$.\epf

A graph $\Gamma$ with diameter at least 2 is called a {\em Terwilliger graph} if for any two vertices $x,y$ with
$d(x,y)=2$, the set $\Gamma(x)\cap \Gamma(y)$ induces a clique. In particular, any distance-regular graphs with
intersection number $c_2=1$ are distance-regular Terwilliger graphs.\\

%In \cite{terw-85}, P. Terwilliger has shown the following result.

%\begin{theorem}(cf. \cite[Theorem 1.16.3]{bcn})\label{T-thm}{\ \\}
%Let $\Gamma$ be an amply regular Terwilliger graph with diameter $D\geq 2$ and parameters $(k,\lambda,\mu)$ where $\mu\ge 2$. Then for any vertex $x\in V(\Gamma)$, the local graph of $x$ is an $\alpha$-clique extension of a strongly regular Terwilliger graph for some integer $\alpha\ge 1$.
%\end{theorem}

\begin{theorem}\label{terwthm}
Fix an integer $m\geq 2$. Then there are only finitely many
distance-regular Terwilliger graphs with smallest eigenvalue at least $-m$ and intersection number
$c_2\ge 2$.
\end{theorem}
\pf For given integer $m\geq 2$, suppose that $\Gamma$ is a
distance-regular Terwilliger graph with $c_2\geq 2$, valency $k$
and  smallest eigenvalue $\theta \geq -m$. By \cite[Theorem
1.16.3]{bcn},  for any vertex $x\in V(\Gamma)$, the local graph
$\Sigma$ of $x$ is the $\alpha$-clique extension of a non-complete
connected strongly regular Terwilliger graph $\Delta$ for an
integer $\alpha \geq 1$. Let $\eta$ be the smallest eigenvalue of
$\Delta$. Then by Proposition \ref{a-clique-ev}, $\alpha(\eta
+1)-1$ is an eigenvalue of $\Sigma$ and hence $\alpha(\eta+1) -1
\geq -m$ as $\Sigma$ is a subgraph of $\Gamma$ and $\Gamma$ has
smallest eigenvalue at least $-m$. As every connected non-complete
strongly regular graph has smallest eigenvalue at most
$\frac{-1-\sqrt{5}}{2}$, we obtain
\[\frac{1-m-\alpha}{\alpha}\le \eta\leq \frac{-1-\sqrt{5}}{2},\]
which immediately gives us
\begin{equation}\label{terw-ev-bd}
\alpha\leq \frac{2(m-1)}{\sqrt{5}-1}.
\end{equation}
By Lemma \ref{geometricpar}, any geometric strongly regular graph
satisfies $c_2 \geq 2$ and thus it is not a
 Terwilliger graph by Lemma \ref{tau} (ii). Hence, it follows by Theorem \ref{neum} that there are only finitely many
 non-complete strongly regular Terwilliger graphs with smallest eigenvalue at least $
 -m$. Hence, the valency $k$ is bounded
above by a function of $m$ by $k = |V(\Sigma)| = \alpha
|V(\Delta)|$ and (\ref{terw-ev-bd}). The theorem now immediately
follows from Theorem \ref{ivanov}. \epf

\noindent {\em Proof of Theorem \ref{main}:} For given integer
$m\geq 2$, suppose that $\Gamma$ is a coconnected non-geometric
distance-regular graph with diameter $D\ge 2$, smallest eigenvalue
at least $-m$ and intersection number $c_2$ satisfying $c_2 \ge
2$. If $\Gamma$ is a Terwilliger graph, then we are done by
Theorem \ref{terwthm}. Now we may assume that $\Gamma$ contains an
induced quadrangle. Proposition \ref{diam-m} implies that $D <
m^2$, from which the theorem follows by Theorem \ref{main-diameter
given}. \epf

\section{Geometric Distance-Regular Graphs}

In this section, we will discuss geometric distance-regular graphs
in more detail.

Let us first return to the strongly regular graph case.

Recall that a {\em partial geometry of order $(s,t, \alpha)$} is a
partial geometry $\mathcal{G}= (\mathcal{P}, \mathcal{L})$ of
order $(s,t)$ such that for each $p \in \mathcal{P}$ and $L \in
\mathcal{L}$ satisfying $p\not \in L$, there are exactly $\alpha$
lines through $p$ that meet $L$. For any partial geometry
$\mathcal{G}=(\mathcal{P}, \mathcal{L})$ of order $(s,t,\alpha)$
with $t \geq 1$, the point graph of $\mathcal{G}$ (i.e its vertex
set is $\mathcal{P}$, and two vertices are adjacent if they lie on
a common line) is a strongly regular graph with intersection
numbers $k= s(t+1), b_1 = (s-\alpha +1)t, c_2 = \alpha(t+1)$ and
distinct eigenvalues $k,s-1, -t-1$, see for example
\cite[Additional result 20f on p.162]{biggs}. As Delsarte cliques
in these graphs have size $s+1$, it follows that each line (each
line is considered as the set of points which lie on it) forms a
Delsarte clique and hence the point graph of a partial geometry of
order $(s, t, \alpha)$ with $t\geq 1$ is geometric. On the other
hand, it is easy to see that a geometric distance-regular graph of
diameter two, is the point graph of a partial geometry of order
$(s,t, \alpha)$ for some $s, t, \alpha$. So these two notions are
equivalent for distance-regular graphs of diameter two.

As we already mentioned in Section \ref{intro}, A. Neumaier showed
that for fixed integer $m \geq 2$, except for a finite number of
graphs, all geometric strongly regular graphs with a given
smallest eigenvalue are either Latin square graphs or Steiner
graphs. Now the Latin square graphs (the Steiner graphs,
respectively) arrive as the point graph of a partial geometry of
order $(s,t,\alpha)$ with $\alpha$ equals to $t+1$ ($t$,
respectively), see for example \cite[Section 4]{cameronsrg}.

As it follows from the above discussion that the point graph of a
partial geometry of order $(s, t, \alpha)$ is geometric with
$\psi_1 = \alpha, \tau_2 = t+1$ and smallest eigenvalue
$-t-1$, A. Neumaier's result implies that geometric strongly
regular graphs with fixed smallest eigenvalue $-t-1 \leq -2$
satisfy $\psi_1 \in \{\tau_2, \tau_2 -1\}$, except for a finite
number of cases.

Now, we return to geometric distance-regular graphs with diameter
at least three. In Lemma \ref{tau} (i), we have shown that
$c_2\geq 2$ implies $\psi_1 \leq \tau_2$. For $c_2 =1$, we have
$\psi_1 = \tau_2 = 1$.

In the following result, we characterize geometric
distance-regular graphs with diameter at least three satisfying
$\psi_1=\tau_2 \geq 2$.  Note that the situation here is completely
different from the case of geometric strongly regular where such a
characterization is impossible.

\begin{theorem}\label{grassmann}
Fix an integer $m\geq 2$. Suppose that $\Gamma$ is a geometric
distance-regular graph with smallest eigenvalue $-m$ and diameter at least three. If $\psi_1 = \tau_2\geq 2$, then one of the following holds.\\
(i) $\psi_1= 2$ and $\Gamma$ is a Johnson graph.\\
(ii) $\psi_1= 2$ and $\Gamma$ is the folded Johnson graph $\overline{J}(4D,2D)$ where $D\geq 3$.\\
(iii) $\psi_1 \geq 3$ and $\Gamma$ is a  Grassmann graph defined over the field $\mathbb{F}_{\psi_1-1}$. \\
(iv) $\psi_1 \geq 3$ and $k < \psi_1 (\psi_1 - 1) m < m^3$.
\end{theorem}
\pf (i)-(ii): Assume that $\tau_2 = \psi_1 = 2$. Then $c_2 = 4$
and the subgraph induced by the common neigbours of two vertices
at distance 2 forms a quadrangle. Now, by \cite[Theorem
9.1.3]{bcn}, the graph $\Gamma$ is either a Johnson graph or a
quotient $J(2s,s)/ \Pi$ for some integer $s\geq 1$, where $\Pi$ is
a uniformly regular partition of $J(2s,s)$ and each part of $\Pi$
has size $2$. By \cite[Theorem 11.1.6]{bcn}, the partition $\Pi$
must be a completely regular partition of $J(2s,s)$, and hence
each part is a completely regular code of size 2 in $J(2s,s)$. If
$C=\{x,y\}$ is a completely regular code of size $2$, then either
$d(x,y)=1$ or $d(x,y)=s$ holds. The case $d(x,y)=1$ only occurs if
either $s=1$ or $a_1=0$. Hence the only completely regular codes
of size 2 in the Johnson graph $J(2s,s)$ are the antipodal pairs
(i.e., $d(x,y)=s$), and this shows that if $\Gamma$ is not a
Johnson graph then it has to be a folded Johnson graph. For a
folded Johnson graph $\overline{J}(2s,s)$, it is geometric with
$c_2 =4$ if
and only if $s$ is even and $s \geq 6$ (cf. \cite[Section 3]{DCG}). This proves (i)-(ii).\\
(iii)-(iv): Suppose that $\tau_2 =\psi_1 \geq 3$ holds. Then by applying \cite[Theorem 9.3.9]{bcn}
(a result of D. K. Ray-Chaudhuri and A. P. Sprague \cite{ray}) with $q:=m-1$, we find that $\Gamma$ is a Grassmann graph
or $|C|  \leq (m-1)^2 + m-1$ for any Delsarte clique $C$. This shows that one of  (iii) and (iv) holds.
\epf

The following corollary follows from Proposition \ref{geom-diam} and Theorem \ref{grassmann}

\begin{cor}\label{corgras}
Fix an integer $m\geq 2$. Suppose that $\Gamma$ is a geometric
distance-regular graph with smallest eigenvalue $-m$ and
diameter $D\geq 3$. If $\psi_1 = \tau_2 \geq 2$ holds, then the
graph $\Gamma$ is either a Johnson graph, a folded Johnson graph,
a Grassmann graph or the number of vertices is bounded above by a
function of $m$.
\end{cor}

Note that the Hamming graph $H(e,D)$ and the bilinear forms graph
$H_q(n,D)$ are geometric distance-regular graphs with $\psi_1 = 1
= \tau_2 -1$ and $\psi_1 = q = \tau_2-1$, respectively. In
\cite{egawa}, Y. Egawa showed that the Hamming graphs are
characterized by its intersection numbers, and K. Metsch
\cite{metsch-99} showed that the bilinear forms graph $H_q(n,D)$
is characterized by its intersection numbers if $n \geq D+4 \geq 7$, by
generalizing results of T. Y. Huang \cite{huang-87} and H. Cuypers
\cite{cuypers}. As far as the authors know, there is no local
geometric characterization known
(as in the case of the Johnson and Grassmann graphs) of neither the bilinear forms graphs nor the Hamming graphs.\\

We close this section with three conjectures concerning geometric
distance-regular graphs.

\begin{conj} \label{conj1}
For fixed integer $m\geq 2$, there are only finitely many coconnected geometric distance-regular graphs
with smallest eigenvalue $-m$ and $\psi_1 \leq  \tau_2-2$.
\end{conj}

\begin{conj}\label{conj2}
For a fixed integer $m\geq 2$, any geometric distance-regular graph
with smallest eigenvalue $-m$, diameter $D\geq 3$ and $c_2\geq
2$ is either a Johnson graph, a Grassmann graph, a Hamming graph,
a bilinear forms graph or the number of vertices is bounded above
by a function of $m$.
\end{conj}

\begin{conj}\label{conj3}
For a fixed integer $m\geq 2$, the diameter of a geometric distance-regular graph with smallest eigenvalue $-m$ and valency at least three
is bounded above by a function of $m$.
\end{conj}

\begin{remark}
(i) Conjecture \ref{conj2} is shown for $\psi_1= \tau_2 \geq 2$ in
Corollary \ref{corgras}. It is known that it is true for $m=2$, but in the case $c_2 =1$ we also have the even polygons. This follows from the result of
Cameron et al. for regular graphs with the smallest eigenvalue $-2$,
as mentioned in the introduction, and the classification of distance-regular line-graphs, see \cite[Theorem 4.2.16]{bcn}. The
second author \cite{bang} has shown it for $m=-3$.\\
(ii) Conjecture \ref{conj3} is shown for $c_2 \geq 2$ in
Proposition \ref{geom-diam}.  In \cite[Problem
3.1.1]{suzukinotes}, H. Suzuki asks a similar question.
\end{remark}

\begin{center}
{\bf Acknowledgements}
\end{center}
The first author was partially supported by a grant of the Korea
Research Foundation funded by the Korean Government (MOEHRD) under
grant number KRF-2008-314-C00007. The second author was supported
by the Korea Research Foundation Grant funded by the Korean
Government(MOEHRD, Basic Research Promotion Fund)
KRF-2008-359-C00002. This support is greatly appreciated. We would
like to thank the anonymous referees for their comments as their
comments greatly improved the paper. Also we would like to thank Jong Yook Park for his careful reading of the paper.


\vspace{5mm}


\begin{thebibliography}{99}

\bibitem{bang}
S.~Bang, Distance-regular graphs with smallest eigenvalue at least
$-3$ and $c_2\geq 2$, in preparation.

\bibitem{DCG}
S.~Bang, A.~Hiraki and J.~H.~Koolen, Delsarte clique graphs, \textit{European J. Combin.} {\bf 28} (2007), 501--516.

\bibitem{banitoconj}
S.~Bang, J.~H.~Koolen and V.~Moulton, There are only finitely many distance-regular graphs of given valency greater than two, in preparation.
%\bibitem{banito}
%E.~Bannai and T.~Ito,
%Algebraic combinatorics I: Association schemes,
%\textit{Benjamin/Cummings, Menlo Park, CA}, 1984.

\bibitem{biggs}
N.~Biggs, Algebraic Graph Theory, Second edition, \textit{Cambridge University Press, Cambridge}, 1993.

\bibitem{bose}
R.~C.~Bose, Strongly regular graphs, partial geometries and partially balanced designs, \textit{Pacific J. Math.} {\bf 13} (1963), 389--419.

\bibitem{bcn}
A.~E.~Brouwer, A.~M.~Cohen and A.~Neumaier,
Distance-regular graphs, \textit{Springer-Verlag, Berlin,} 1989.

\bibitem{cameronsrg}
P.~J.~Cameron, Strongly regular graphs, in: Topics in algebraic graph theory, L.W. Beineke and R.J. Wilson eds., Encyclopedia of Mathematics and Its Applications {\bf 102}, Cambridge University Press, Cambridge, 2004, 203--221.

\bibitem{cameron et al}
P.~J.~Cameron, J.~M.~Goethals, J.~J.~Seidel and E.~E.~Shult, Line graphs, root systems, and elliptic geometry, \textit{J. Algebra} {\bf 43} (1976), 305--327.

%\bibitem{cuypersgrass}
%H.~Cuypers, The dual of Pasch's axiom, \textit{European J. Combin.} {\bf 13} (1992), 15--31.

\bibitem{cuypers}
H.~Cuypers, Two remarks on Huang's characterization of the bilinear forms graphs, \textit{European J. Combin.} {\bf 13} (1992), 33--37.

\bibitem{clerck-06}
F.~De~Clerck, S.~De~Winter, E.~Kuijken and C.~Tonesi,
Distance-regular $(0,\alpha)$-reguli, \textit{Des. Codes Cryptogr.} {\bf 38} (2006), 179--194.

\bibitem{Delsarte}
P.~Delsarte, An algebraic approach to the association schemes of coding theory, \textit{Philips Res. Rep. Suppl.} {\bf 10} (1973).

\bibitem{egawa}
Y.~Egawa, Characterization of $H(n,q)$ by the parameters, \textit{J. Combin. Theory Ser. A} {\bf 31} (1981), 108--125.

%\bibitem{fdf}
%D.~G.~Fon-Der-Flaass, New prolific constructions of strongly regular graphs, \textit{Adv. Geom.} {\bf 2} (2002), 301--306.

\bibitem{godsil-combin}
C.~D.~Godsil, Algebraic combinatorics, Chapman and Hall Mathematics Series, \textit{Chapman and Hall, New York,} 1993.

\bibitem{godsil-93}
C.~D.~Godsil, Geometric distance-regular covers, \textit{New Zealand J. Math.} {\bf 22} (1993), 31--38.

\bibitem{huang-87}
T.~Y.~Huang, A characterization of the association schemes of
bilinear forms, \textit{European J. Combin.} {\bf 8} (1987),
159--173.

\bibitem{aaivanov}
A.~A.~Ivanov,
Bounding the diameter of a distance-regular graph,
\textit{Soviet Math. Dokl.} {\bf 28} (1983), 149--152.

\bibitem{metsch-99}
K.~Metsch, On a characterization of bilinear forms graphs, \textit{European J. Combin.} {\bf 20} (1999), 293--306.

\bibitem{neumaier-79}
A.~Neumaier, Strongly regular graphs with smallest eigenvalue $-m$, \textit{Arch. Math. (Basel)} {\bf 33} (1979/1980), 392--400.

\bibitem{ray}
D.~K.~Ray-Chaudhuri and A.~P.~Sprague, A combinatorial characterization of attentuated spaces, \textit{Utilitias Math.} {\bf 15} (1979), 3--29.

\bibitem{suzuki-semiregular}
H.~Suzuki, Distance-semiregular graphs, \textit{
Algebra Colloq.} {\bf 2} (1995), 315--328.

\bibitem{suzukinotes}
H.~Suzuki, An introduction to distance-regular graphs, p. 57-132 in: Three lectures in Algebra, Sophia Kokyuroku in Mathematics 41 (K. Shinoda ed.) (1999).

\bibitem{terw-85}
P.~Terwilliger, Distance-regular graphs with girth 3 or 4, I, \textit{J. Combin. Theory Ser. B} {\bf 39} (1985), 265--281.

\bibitem{wilson-74}
R.~M.~Wilson, Nonisomorphic Steiner triple systems, \textit{Math. Z.} {\bf 135} (1973/74), 303--313.

\end{thebibliography}
\end{document}